\documentclass{amsart}
\input epsf

\usepackage{fullpage}

\usepackage{amsfonts,amsthm,amsmath,amssymb,latexsym}

\def\H{\mathcal{H}}

\def\S{\mathcal{H}}

\begin{document}
\title{Extremal maps of the universal hyperbolic solenoid}
\author{Adam Epstein, Vladimir Markovic and Dragomir \v{S}ari\'{c}}

\address{Institute of Mathematics, University of Warwick, CV4 7AL Coventry, UK}
\email{adame@maths.warwick.ac.uk}

\address{Institute of Mathematics, University of Warwick, CV4 7AL Coventry, UK}
\email{markovic@maths.warwick.ac.uk}

\address{Institute for Mathematical Sciences, Stony Brook University,
Stony Brook, NY 11794-3660} \email{saric@math.sunysb.edu}


\date{\today}

\begin{abstract}

We show that the set of points in the Teichm\"uller space of the
universal hyperbolic solenoid which do not have a Teichm\"uller
extremal representative is generic (that is, its complement is the set of the first
kind in the sense of Baire). This is in sharp contrast with the Teichm\"uller space of a
Riemann surface where at least an open, dense subset has
Teichm\"uller extremal representatives. In addition, we provide a
sufficient criteria for the existence of Teichm\"uller extremal
representatives in the given homotopy class. These results
indicate that there is an interesting theory of extremal (and
uniquely extremal) quasiconformal mappings on hyperbolic
solenoids.

\end{abstract}

\maketitle

%VT
\thispagestyle{empty}
\def\IMSmarkvadjust{0 pt}
\def\IMSmarkhadjust{0 pt}
\def\IMSmarkhpadding{0 pt}
\def\IMSpubltext{Published in modified form:}
\def\SBIMSMark#1#2#3{
 \font\SBF=cmss10 at 10 true pt
 \font\SBI=cmssi10 at 10 true pt
 \setbox0=\hbox{\SBF \hbox to \IMSmarkhpadding{\relax}
                Stony Brook IMS Preprint \##1}
 \setbox2=\hbox to \wd0{\hfil \SBI #2}
 \setbox4=\hbox to \wd0{\hfil \SBI #3}
 \setbox6=\hbox to \wd0{\hss
             \vbox{\hsize=\wd0 \parskip=0pt \baselineskip=10 true pt
                   \copy0 \break%
                   \copy2 \break% 
                   \copy4 \break}}
 \dimen0=\ht6   \advance\dimen0 by \vsize \advance\dimen0 by 8 true pt
                \advance\dimen0 by -\pagetotal
	        \advance\dimen0 by \IMSmarkvadjust
 \dimen2=\hsize \advance\dimen2 by .25 true in
	        \advance\dimen2 by \IMSmarkhadjust

%
%   Check for publication info
%
%  \newread\jref
  \openin2=publishd.tex
  \ifeof2\setbox0=\hbox to 0pt{}
  \else 
     \setbox0=\hbox to 3.1 true in{
                \vbox to \ht6{\hsize=3 true in \parskip=0pt  \noindent  
                {\SBI \IMSpubltext}\hfil\break
                \input publishd.tex 
                \vfill}}
  \fi
  \closein2
  \ht0=0pt \dp0=0pt
 \ht6=0pt \dp6=0pt
 \setbox8=\vbox to \dimen0{\vfill \hbox to \dimen2{\copy0 \hss \copy6}}
 \ht8=0pt \dp8=0pt \wd8=0pt
 \copy8
 \message{*** Stony Brook IMS Preprint #1, #2. #3 ***}
}

\SBIMSMark{2006/02}{April 2006}{}

\section{Introduction}

The Teichm\"uller space $T(S)$ of a Riemann surface $S$ consists
of all marked complex structures on $S$. A {\it marked complex
structure} on $S$ is a homotopy class of quasiconformal maps from
$S$ to an arbitrary Riemann surface up to post composition by a
conformal map. The homotopy class $[id]$ of the identity map
$id:S\to S$ is the {\it basepoint} of $T(S)$. The distance between
the basepoint and the homotopy class $[f]$ of a quasiconformal map
$f:S\to S_1$ is the infimum of the logarithms of quasiconformal
constants over all quasiconformal maps in the marked complex
structure($\equiv$homotopy class $[f]$) of $f$. A quasiconformal
map $f_1:S\to S_1$ in the homotopy class $[f]$ is called {\it
extremal} if its quasiconformal constant $K(f_1)$ is equal to the
infimum of the quasiconformal constants over all maps in $[f]$. In
this case, the (Teichm\"uller) distance between $[id]$ and $[f]$
is simply:
$$
dist([id],[f])=\log K(f_1).
$$

\vskip .3 cm

\footnotetext[1]{The third author was partially supported by NSF
grant DMS-0505652.}

If $\mu$ is a Beltrami coefficient of an extremal map $f:S\to
S_1$, then each quasiconformal map $f^{t\mu}$ with the Beltrami
coefficient $t\mu$, $\frac{-1}{\|\mu\|_{\infty}}
<t<\frac{1}{\|\mu\|_{\infty}}$, is extremal as well. The path
$t\mapsto [f^{t\frac{|\varphi |}{\varphi}}]\in T(S)$ is a
geodesic for the above distance.

\vskip .3 cm

Teichm\"uller's fundamental result \cite{T} states that each
marked complex structure of an analytically finite (closed with at
most finitely many points deleted) Riemann surface $S$ contains a
unique extremal map with Beltrami coefficient $k\frac{|\varphi
|}{\varphi}$, where $\varphi$ is a holomorphic quadratic
differential on $S$ and $0<k<1$. The natural parameter for
$\varphi$ partitions $S$ into Euclidean rectangles and the
extremal map is an affine stretching on each rectangle.
Teichm\"uller's theorem is a highly non-trivial generalization of
a result of Gr\"otzsch concerning a single rectangle. We say that
such extremal maps, their Beltrami coefficients and their
corresponding geodesics $t\mapsto [f^{t\frac{|\varphi
|}{\varphi}}]$ are of {\it Teichm\"uller type}.

\vskip .3 cm

The Teichm\"uller theorem completely answers all questions about
extremal and uniquely extremal quasiconformal mappings for
analytically finite Riemann surfaces. However, for Riemann
surfaces that are not analytically finite there exists a rich
theory about extremal mappings. A modern approach to extremal maps
for arbitrary Riemann surfaces started with Reich and Strebel.
They showed that Teichm\"uller maps are extremal on arbitrary
Riemann surfaces by generalizing the original approach of
Gr\"otzsch to Riemann surfaces. Using, what is today called
Reich-Strebel inequality and results of Hamilton \cite{H} and
Krushkal \cite{K}, they characterized extremal quasiconformal maps
\cite{RS}. It is worth noting that every marked complex structure
of a Riemann surface contains an extremal map by the
pre-compactness of a family of normalized $K$-quasiconformal maps.
Strebel \cite{Str} showed that not every extremal map is of
Teichm\"uller type and that there could be more than one extremal
map in a given homotopy class on analytically infinite Riemann
surfaces. The corresponding characterization for uniquely extremal
maps has been obtained by Bozin, Lakic, Markovic and Mateljevic
\cite{BLMM} (see \cite{Mar}, \cite{Re}, \cite{Bi}, \cite{Y} for
some applications of these results).

\vskip .3 cm

Sullivan \cite{Sul} introduced the universal hyperbolic solenoid
$\S$ as the inverse limit of the system of unbranched finite
degree covers of a compact surface. The universal hyperbolic
solenoid $\S$ is a compact space which is locally homeomorphic to
a 2-disk times a Cantor set. Sullivan \cite{Sul} introduced a
complex structure on $\S$ and showed that the Teichm\"uller space
$T(\S )$ of the solenoid is a separable complex Banach manifold.
Nag and Sullivan \cite{NS} observed that $T(\S )$ embeds in the
universal Teichm\"uller space as a closure of the union of
Teichm\"uller spaces of all compact Riemann surfaces. In the
remark at the end of Section 2, we use the embedding to give an
alternative definition of $T(\S )$ as a subset of the universal
Teichm\"uller space and the metric on it is given in terms of this
subset (this description uses terms more familiar with the standard
Teichm\"uller theory). The questions that we consider can be directly restated in
these terms.) The Teichm\"uller distance on $T(\S )$ is defined
similar to above and we consider the questions about the existence
and the structure of extremal quasiconformal maps. Unlike for
Riemann surfaces, the existence of extremal maps on the solenoid
is not guaranteed, and in fact it is an interesting open problem
whether they always exist. In \cite{Sa}, it is showed that
Teichm\"uller type maps are uniquely extremal. The question was
raised whether each marked complex structure contains a
Teichm\"uller map. If the answer were yes, this would show both
the existence and the uniqueness of the extremal maps.

\vskip .3 cm

There are reasons why one would expect to have a positive answer.
The solenoid $\S$ is a compact space and $T(\S )$ is the closure
of the union of Teichm\"uller spaces of all compact Riemann
surfaces. (Recall that by Teichm\"uller's theorem each marked
complex structure on a compact surface contains a Teichm\"uller
map.) On the other hand, Lakic \cite{Lak} showed that the set of
points in the Teichm\"uller space of an analytically infinite
Riemann surface which contain a Teichm\"uller extremal map is open
and dense. The set of points which do not have a Teichm\"uller
representative is nowhere dense and closed, in particular of first
kind in the sense of Baire. By analogy, one would expect the
subset of $T(\S )$ without Teichm\"uller extremal maps is at most
of first kind, if not empty. We obtain a somewhat unexpected result:

\vskip .4 cm

\paragraph{\bf Theorem 1} {\it The set of points in the Teichm\"uller space
$T(\S )$ of the universal hyperbolic solenoid $\S$ which do not
have a Teichm\"uller extremal representative is generic in
$T(\S)$. That is, the set of points that do have a Teichm\"uller
representative is of the first kind in the sense of Baire with
respect to the Teichm\"uller metric.}

\vskip .4 cm

By the above theorem, a large set of points in $T(\S )$ do not
have a Teichm\"uller extremal representative. It is interesting to
determine when a given marked complex structure has a
Teichm\"uller extremal representative. A sufficient condition for
the case of infinite Riemann surfaces is given by Strebel's Frame
Mapping Condition \cite{Str1}. This condition depends on the
existence of ends of the Riemann surface. Since the universal
hyperbolic solenoid $\S$ is a compact space, it has no ends.
However, we obtain a sufficient condition, where a given complex structure has a Teichm\"uller representative if it
is well approximated with locally constant (rational) complex structures (see below for the definition of locally constant complex structures).

\vskip .3 cm

A complex structure on an arbitrary compact Riemann surface lifts
to a special complex structure on the solenoid $\S$, namely a
transversely locally constant complex structure. Each transversely
locally constant marked complex structure on $\S$ has a
Teichm\"uller extremal representative which comes by lifting the
extremal Teichm\"uller representative from the surface. Therefore,
it is interesting to consider only non transversely locally
constant complex structures. In Section 3 we define a notion of
being {\it well-approximated by transversely locally constant
complex structures}. (Recall that the transversely locally
constant complex structures are dense in $T(\S )$ by \cite{Sul}.
Therefore each non transversely locally constant complex structure
is approximated by transversely locally constant complex
structures.) We show:

\vskip .4 cm

\paragraph{\bf Theorem 2} {\it If a non locally transversely constant
marked complex structure is well-approximated by transversely
locally constant complex structures then it contains a
Teichm\"uller extremal representative.}

\vskip .4 cm

The immediate consequence of Theorem 1 and Theorem 2 is

\vskip .4 cm

\paragraph{\bf Corollary 1} {\it The set of points in $T(\S )$
which are not well-approximated by transversely locally constant
marked complex structures is generic in $T(\S )$.}

\vskip .4 cm

This is analogous with the fact that the set of real numbers
which are not well approximated by rational numbers is of full
measure. We remark that the existence of extremal maps in
arbitrary marked complex structure is open. Further, the existence
of a geodesic connecting the basepoint with any marked complex
structure is implied by the existence of an extremal map but it is
not necessarily equivalent to it. This is illustrated by an
example of L. Zhong \cite{Zh} for Riemann surfaces.

\section{Preliminaries}

We define the {\it universal hyperbolic solenoid} $\H$ introduced
by Sullivan \cite{Sul}, see also \cite{Odd} and \cite{Sa}. Let
$(S,x)$ be a fixed compact surface of genus at least $2$ with a
basepoint. Consider all finite sheeted unbranched coverings
$\pi_i:(S_i,x_i)\to (S,x)$ such that $\pi_i(x_i)=x$. There is a
natural {\it partial ordering} $\geq$ given by $(S_j,x_j)\geq
(S_i,x_i)$ whenever there exists a finite unbranched covering
$\pi_{i,j}:(S_j,x_j)\to (S_i,x_i)$ such that $\pi_{i,j}(x_j)=x_i$
and $\pi_j=\pi_{i}\circ\pi_{i,j}$. Given any two finite unbranched
covers $(S_i,x_i)$ and $(S_j,x_j)$ of $(S,x)$, there exists a
third finite unbranched cover $(S_k,x_k)$ which covers both such
that $\pi_i\circ\pi_{i,k}=\pi_j\circ\pi_{j,k}=\pi_k$,
$\pi_{i,k}(x_k)=x_i$ and $\pi_{j,k}(x_k)=x_j$. In other words,
$(S_k,x_k)\geq (S_i,x_i)$ and $(S_k,x_k)\geq (S_j,x_j)$. Thus the
system of covers is inverse directed and the inverse limit is well
defined. The universal hyperbolic solenoid is by definition
$$
\S =\lim_{\leftarrow} (S_i,x_i).
$$

\vskip .3 cm

We give an alternative definition for the universal hyperbolic
solenoid $\S$. Denote by $\Delta$ the unit disk. Let $G$ be a
Fuchsian group such that $\Delta /G$ is a compact Riemann surface
of genus at least two. Let $G_n$ be the intersection of all
subgroups of $G$ of index at most $n$. Then $G_n$ is a finite
index characteristic subgroup of $G$. The {\it profinite metric}
on $G$ is defined by
$$
d(\alpha ,\beta )=\max_{\alpha\beta^{-1}\in G_n}\frac{1}{n}
$$
for all $\alpha ,\beta\in G$. The profinite completion $\hat{G}$
of $G$ is a compact topological group homeomorphic to a Cantor
set. The action of $G$ on the product $\Delta\times\hat{G}$ is
defined by $\gamma (z,t):=(\gamma (z),t\gamma^{-1})$ for all
$\gamma\in G$ and $(z,t)\in \Delta\times\hat{G}$. The universal
hyperbolic solenoid is $\S :=(\Delta\times\hat{G})/G$. For more
details on this definition see \cite{Odd}.

\vskip .3 cm

The solenoid $\S$ is a compact space locally homeomorphic to a
2-disk times a Cantor set. Path components of $\S$ are called {\it
leaves}. Each leaf is dense in $\S$ and homeomorphic to the unit
disk. The profinite group completion $\hat{G}$ supports a unique
left and right translation invariant measure $m$ of full support
(the Haar measure). The Haar measure $m$ induces a holonomy
invariant measure on the solenoid $\S =(\Delta\times\hat{G})/G$.
This measure allows for the integration of quantities which induce
local measures on leaves of $\S$, e.g. the absolute value of a
quadratic differential.

\vskip .3 cm

The complex structure on $\S$ is given by an assignment of
holomorphic charts on leaves(making the leaves holomorphic to the
unit disk) which vary continuously for the transverse variation in
local charts \cite{Sul}. From the results of Candel \cite{Can}, it
follows that $\S$ supports hyperbolic metric for each conformal
class (induced by a complex structure). Note that by fixing a
Fuchsian group $G$ the solenoid $\S=(\Delta\times\hat{G})/G$ has
already induced complex structure and hyperbolic metric from the
unit disk $\Delta$ (see \cite{Odd}). The induced complex structure
is locally constant in the transverse direction and all locally
transversely constant complex structures on $\S$ arise in this way
\cite{Sul}, \cite{NS}.

\vskip .3 cm

The Teichm\"uller space $T(\S )$ consists of all differentiable
quasiconformal maps (which are continuous for the transverse
variation in local charts) from fixed complex solenoid $\S
=(\Delta\times\hat{G})/G$ onto an arbitrary complex solenoid
modulo homotopy and post-composition by conformal maps. The
requirement of the continuity for the transverse variation in
local charts can be achieved by requiring the differentiable maps
to vary continuously in the $C^{\infty}$-topology on $C^{\infty}$
maps (for more details see \cite{Sa}). Equivalently, the
Teichm\"uller space $T(\S )$ is the space of all smooth Beltrami
coefficients which are continuous for the transverse variations in
the $C^{\infty}$-topology modulo the above condition. Thus a point
in $T(\S )$ is an equivalence class $[\mu ]$ of a smooth Beltrami
coefficient $\mu$ on $\S$. The main point is that the leafwise
Beltrami equations give a transversely continuous solution (see
\cite{Sa}). It is also possible to weaken the condition on
differentiability of the Beltrami coefficients as long as they are
leafwise equivalent(as elements of the universal Teichm\"uller
space) to the restriction on leaves of a smooth Beltrami
coefficient on $\S$ (for more details see \cite{Sa}). A Beltrami
coefficient $\mu$ on $\S$ is {\it extremal} if
$\|\mu\|_{\infty}=\inf_{\nu\in [\mu ]} \|\nu\|_{\infty}$. A
Beltrami coefficient $\mu$ is of {\it Teichm\"uller} type if $\mu
=k\frac{|\varphi |}{\varphi}$ for a holomorphic quadratic
differential $0\neq\varphi$ and $0<k<1$. They are corresponding to
the similar notions for quasiconformal maps given in Introduction.

\vskip .3 cm

The space $A(\S )$ consists of all holomorphic quadratic
differentials on $\S$ which are continuous for the transverse
variation in local charts. The Bers norm of $\varphi\in A(\S )$ is
given by $\|\varphi\|_{Bers}:=\|\varphi\rho^{-2}\|_{\infty}$,
where $\rho$ is the hyperbolic length element on the leaves of
$\S$. The space $A(\S )$ is a complex Banach space for the Bers
norm and the closure of $A(\S )$ for the $L^1$-norm given by
$\|\varphi\|_{L^1}=\int_{\S}|\varphi |dm$ is the space of all
integrable quadratic differentials which are holomorphic(and
defined) on almost all leaves of $\S$ without transverse
continuity requirement(see \cite{Sa}).

\vskip .3 cm

The space $N(\S )$ of {\it infinitesimally trivial Beltrami
differentials} consists of all smooth Beltrami differentials $\mu$
such that $\int_{\S}\mu\varphi dm=0$, for all $\varphi\in A(\S )$.
In fact, a smooth Beltrami differential $\mu$ is infinitesimally
trivial if and only if there exists a path of smooth Beltrami
coefficients $t\mapsto\nu_t$ such that $\nu_t=t\mu+o(t)$ and
$\nu_t$ is a trivial deformation of $\S$, i.e. the quasiconformal
map $f^{\nu_t}$ is homotopic to the identity (see \cite{Sa}).

\vskip .3 cm

The tangent space to the Teichm\"uller space $T(\S )$ at the base
point is given by the space $L^{\infty}(\S )$ of smooth Beltrami
differentials on $\S =(\Delta\times\hat{G})/G$ modulo the space
$N(\S )$ of infinitesimally trivial Beltrami differentials(for
details see \cite{Sa}). (Our terminology assumes that any Beltrami
coefficient $\mu$ satisfies $\|\mu\|_{\infty}<1$, while any
Beltrami differential $\mu$ satisfies $\|\mu\|_{\infty}<\infty$.)
A Beltrami differential $\mu$ on $\S$ is {\it infinitesimally
extremal} if $\|\mu\|_{\infty}=\inf_{\nu}\|\nu\|_{\infty}$, where
the infimum is over all $\nu$ such that $\mu -\nu\in N(\S )$. The
tangent space $L^{\infty}(\S )/N(\S )$ is Banach in the quotient
topology even though $L^{\infty}(\S )$ and $N(\S )$ are not
complete (see \cite{Sa}). There is a natural pairing between
$L^{\infty}(\S )$ and $A(\S )$ given by
$$
(\mu ,\varphi)\mapsto\int_{\S}\mu\varphi dm.
$$
The pairing descends to the pairing of $L^{\infty}(\S )/N(\S )$
and $A(\S )$. The tangent space $L^{\infty}(\S )/N(\S )$ embeds in
the dual $A(\S )^{*}$ but it is strictly smaller \cite{Sa}.

\vskip .3 cm

\paragraph{\bf Remark} We give an alternative description of the
Teichm\"uller space $T(\S )$. Fix a Fuchsian group $G$ such that
$\Delta /G$ is a compact Riemann surface of genus at least two. A
quasiconformal map $f:\Delta\to\Delta$ is said to be {\it almost
invariant} with respect to $G$ if $\| Belt(f\circ\gamma\circ
f^{-1})\|_{\infty}\to 0$ as $d(\gamma ,id)\to 0$. In other words,
the Beltrami coefficient of $f$ is very close to be invariant
under the push forward by elements of a finite index subgroup
$G_n$ of $G$ of some large index (i.e. $n$ is large). In
particular, a lift of a quasiconformal map from the Riemann
surface $\Delta /G_n$ for any finite index subgroup $G_n$ of $G$
is almost invariant for $G$. The Teichm\"uller space $T(\S )$ is
isomorphic to a subset of the universal Teichm\"uller space
$T(\Delta )$  consisting of all classes with almost invariant (for
$G$) representatives. The distance $dist([id],[f])$ is given by
the infimum of the logarithm of the quasiconformal constants of
all almost invariant maps homotopic to $f$. The question about
extremal representatives can be considered in this setting as
well. However, it appears that working directly on the solenoid
$\S$ is somewhat better suited for our purposes due to the strong
technical tools developed using the Reich-Strebel inequality for
the solenoid \cite{Sa}. In fact, one would presumably be able to
replace integration in the transverse direction by the limit of
the average of integrals over fundamental regions for $G_n$ as
$n\to\infty$.

\section{A sufficient condition for Teichm\"uller maps}

A Teichm\"uller Beltrami coefficient $\mu =k\frac{|\varphi
|}{\varphi}$, for some $0<k<1$ and for some holomorphic quadratic
differential $\varphi\neq 0$ on the solenoid $\H$, is uniquely
extremal in its Teichm\"uller class and it determines a geodesic
$t\mapsto [t\frac{|\varphi |}{\varphi}]$, $t\in (-1/k,1/k)$.
Moreover, this is unique geodesic connecting the basepoint $[0]$
with $[\mu =k\frac{|\varphi |}{\varphi}]$ (see \cite{Sa}).

\vskip .3 cm

If $\mu =k\frac{|\varphi |}{\varphi}$, $k\in\mathbb{R}^{+}$, is a
Beltrami differential, then the linear functional
$\Lambda_{\mu}:\psi\mapsto\int_{\H}\mu\psi dm$, for $\psi\in A(\H
)$, achieves its norm on the vector
$\frac{\varphi}{\|\varphi\|_{L^1}}$. In that case,
$\|\mu\|_{\infty} =k$ is equal to the norm of the functional
$\Lambda_{\mu}$ and any other $\nu$ in the infinitesimal class of
$\mu$(i.e. any $\nu$ such that $\mu -\nu\in N(\S )$) satisfies
$\|\nu\|_{\infty}>\|\mu\|_{\infty}$. In other words, $\mu$ is
uniquely infinitesimally extremal.

\vskip .3 cm

It is not, a priori, clear whether each Teichm\"uller (or
infinitesimal) class contains a Teichm\"uller type Beltrami
coefficient. If this is the case, this would certainly be a nice
situation similar to Teichm\"uller spaces of compact surfaces. On
the other hand, on infinite Riemann surfaces there exist
Teichm\"uller (and infinitesimal) classes of Beltrami coefficients
(and differentials) which do not contain a Teichm\"uller type
Beltrami coefficient (and differential). Strebel \cite{Str} gave a
very useful sufficient condition (called the Frame Mapping
Condition) to determine when a given class contains a
Teichm\"uller type representative.

\vskip .3 cm

We find a sufficient condition for a given Beltrami coefficient
$\mu$ on the solenoid $\H$ to be equivalent to a Teichm\"uller
type Beltrami coefficient in both infinitesimal and Teichm\"uller
classes. We point out that the Strebel's Frame Mapping Condition
depends on the non-compactness of the given Riemann surface,
whereas the solenoid is a compact space. Therefore we need a
different approach. If $\mu$ is a transversely locally constant
Beltrami coefficient then it is a lift of a Beltrami coefficient
on a Riemann surface $S_i$ covering $S$ for the base complex
structure on $\H$($\equiv \Delta\times\hat{G}_i /G_i$, where
$S\equiv \Delta /G_i$, $G_i$ a Fuchsian group). Since on $S_i$ any
Beltrami coefficient is equivalent to a Teichm\"uller type
Beltrami coefficient, by lifting the corresponding holomorphic
quadratic differential on $S_i$ to $\H$, we obtain a Teichm\"uller
coefficient equivalent to $\mu$ (either infinitesimally or in
Teichm\"uller sense). Therefore we restrict our attention to non
transversely locally constant Beltrami coefficients on $\S$ and
look for a sufficient condition.

\vskip .3 cm

Let $S_n\equiv\Delta /G_n$ ($S_1=S$) be a sequence of finite
sheeted coverings of $S$ such that $\cap_{n=1}^{\infty}G_n=\{
id\}$ (we assume that $G_{n+1}<G_n$). One can think about $S_n$ as
an approximating sequence for $\H$.

\vskip .3 cm

Let $\mu$ be a non-trivial (in the Teichm\"uller sense) Beltrami
coefficient on $\Delta /G_n$ such that $[\mu ]=[k\frac{|\varphi
|}{\varphi}]$, for $0<k<1$ and $\varphi$ normalized such that
$\|\varphi\|_{Bers}=1$. Let $B_n:T(\Delta /G_n)\to A(\Delta /G_n)$
be the map given by $B_n([\mu])=k\varphi$, where $\mu$, $k$ and
$\varphi$ are as above. The map $B_n$ is continuous for the
Teichm\"uller  metric on $T(\Delta /G_n)$ and the Bers norm on the
unit ball in $A(\Delta /G_n)$ (see \cite{Gar}).

\vskip .3 cm

\paragraph{\bf Definition 3.1}
 Let $\mu$ be a  Beltrami
coefficient on $\H$ not equivalent to a transversely locally
constant Beltrami coefficient. Let $S_n$ be an increasing sequence
of finite coverings of $S$, $\varphi_n$ a sequence of holomorphic
quadratic differentials on $S_n$ and $\tilde{\varphi}_n$ their
lifts to $\H$. Given a sequence $0<k_n<1$, define a sequence of
Beltrami coefficients $\mu_n=k_n\frac{|\varphi_n|}{\varphi_n}$ on
$S_n$ and their lifts
$\tilde{\mu}_n=k_n\frac{|\tilde{\varphi}_n|}{\tilde{\varphi}_n}$
on $\H$. Assume that $[\tilde{\mu}_n]\to [\tilde{\mu} ]$ for some
$0<k_n<1$. The Teichm\"uller class $[\mu ]$ of the Beltrami
coefficient $\mu$ is {\it well-approximated by transversely
locally constant Beltrami coefficients} if there exists a sequence
$\mu_n$ as above such that
$$
\sum_{n=1}^{\infty}\|
B_n([\mu_n])-B_{n+1}([\mu_{n+1}])\|_{Bers}<\infty .
$$

\vskip .3 cm

\paragraph{\bf Proof of Theorem 2} Consider the lifts $\tilde{\varphi}_n$ on
$\H$ of holomorphic quadratic differentials $\varphi_n$ on $S_n$.
Then
$\sum_{n=1}^{\infty}\|\tilde{\varphi}_n-\tilde{\varphi}_{n+1}\|_{Bers}
<\infty$ by the assumption, which implies that $\tilde{\varphi}_n$
converges uniformly to a holomorphic quadratic differential $\psi$
on $\H$. Note that $\psi$ is not a lift of a holomorphic quadratic
differential on $S_n$, for any $n$. Since $\psi
=\tilde{\varphi}_k+\sum_{n=k}^{\infty}
(\tilde{\varphi}_{k+1}-\tilde{\varphi}_k)$ and
$\|\tilde{\varphi}_n\|_{Bers}=1$, we conclude that $\psi\neq 0$.
By the uniform convergence $\tilde{\varphi}_n\to\psi$, we get that
$\mu$ is Teichm\"uller equivalent to $k\frac{|\psi |}{\psi}$, for
$k$ depending on the distance from $[0]$ to $[\mu ]$. $\Box$

\vskip .3 cm

A similar statement can be made for the infinitesimal case. Let
$\mu$ be a Beltrami differential on $\S$ representing a tangent
vector $[\mu ]$ which does not come from lifting a tangent vector
of the Teichm\"uller space of a compact surface, i.e. the coset
$\mu +N(\S )$ does not contain a transversely locally constant
Beltrami differential. Similar to Teichm\"uller classes, it is
also true that transversely locally constant Beltrami
differentials approximate each Beltrami differential on $\S$.
(Recall that each tangent vector is a continuous linear functional
on the space of holomorphic quadratic differentials $A(\S )$ and
the approximation is with respect to the dual norm.)

\vskip .3 cm

Let $S_n$ be a sequence of compact Riemann surfaces
``approximating'' $\S$ as above. Define $B_n':T(S_n)\to A(S_n)$ by
$B_n'([\mu_n ])=k\varphi_n$, where $\mu_n -k\frac{|\varphi_n
|}{\varphi_n}\in N(S_n)$ and $\|\varphi_n\|_{Bers}=1$. We say that
the infinitesimal class of a non transversely locally constant
Beltrami differential $\mu$ on $\S$ is {\it well-approximated with
transversely locally constant Beltrami differentials} if there
exists a sequence of Beltrami differentials $\mu_n$ on $S_n$ whose
lifts $\tilde{\mu}_n$ on $\S$ approximate $\mu$ in the sense of
the linear functionals on $A(\S )$ such that $\sum_{n=1}^{\infty}
\| B_n'([\mu_n])-B_{n+1}'([\mu_{n+1}])\|_{Bers}<\infty$. We obtain
an analogous statement to Theorem 2(and similar proof) for the
infinitesimal class.

\vskip .4 cm

\paragraph{\bf Theorem 2'} {\it If a non transversely locally
constant infinitesimal class of a Beltrami differential on the
universal hyperbolic solenoid is well-approximated by transversely
locally constant infinitesimal classes of Beltrami differentials
then it is infinitesimally equivalent to a Teichm\"uller Beltrami
differential.}

\section{Teichm\"uller classes without Teichm\"uller representatives}

In this section we consider the question of existence of
Teichm\"uller representatives for arbitrary Teichm\"uller classes
in $T(\H )$. We show that Teichm\"uller representative does not
always exist. It is true that there exists a dense subsets of
points in $T(\H )$ which have Teichm\"uller representative by the
density of transversely locally constant structures on $\H$ (which
come from complex structures on finite sheeted covers of $S$).

\vskip .3 cm

Lakic \cite{Lak} showed that even though not all Teichm\"uller
classes of Beltrami coefficients on infinite Riemann surfaces have
Teichm\"uller representative, the one that do form an open, dense
subset of the corresponding Teichm\"uller space. Therefore, for
infinite Riemann surfaces this set is quite large and for finite
Riemann surfaces it equals the whole Teichm\"uller space.

\vskip .3 cm

We show that, quite unexpectedly, for $T(\H )$ the set of points
which do not have Teichm\"uller Beltrami coefficient
representative is generic.  This means that the set of elements in
$T(\H )$  that do have a Teichm\"uller Beltrami representative is contained
in a countable union of closed nowhere dense subsets of $T(\H )$.

\vskip .3 cm

Let $Arg:\mathbb{C}-\{ 0\}\to (-\pi ,\pi]$ be the standard
argument function defined on non-zero complex numbers. Let
$\Delta_r=\{ z\in\mathbb{C};\ |z|<r\}$. Given a holomorphic
function $f:\Delta_1\to\mathbb{C}$ and $0<r< 1$, denote by $\|
f|_{\Delta_r}\|_{Bers}$ the supremum of $|f(z)|\rho^{-2}(z)$ over
$\Delta_r$, where $\rho$ is the hyperbolic length density on the
unit disk $\Delta =\Delta_1$. Let $\| f\|_{Bers}$ denote the Bers
norm, namely the supremum of $|f(z)|\rho^{-2}(z)$ over the unit
disk $\Delta_1$ if it exists. Given a measurable set
$S\subset\Delta_1$, denote by $|S|$ its Euclidean area. In what
follows, we use the following lemma.

\vskip .3 cm

\paragraph{\bf Lemma 4.1}{\it Let $\epsilon >0$, $0<r<1$, $N>M>0$ and
let $\varphi ,\psi$ be two holomorphic functions on the unit disk
$\Delta_1$ such that $\|\varphi \|_{Bers},\|\psi \|_{Bers}\leq N$,
$\|\varphi |_{\Delta_r}\|_{Bers}\geq M$. Assume that there exists
a measurable set $S\subset\Delta_1$ with $|S|=p>0$ such that $\|
Arg(\frac{\psi}{\varphi} )|_S\|_{\infty}\leq\epsilon$. Then there
exists $k>0$ and $\delta (\epsilon ,p,r,M,N)>0$ such that
$$
\| (\psi -k\varphi )|_{\Delta_r}\|_{Bers} <\delta (\epsilon
,p,r,M,N),
$$
where $\delta (\epsilon ,p,r,M,N)\to 0$ as $\epsilon\to 0$ for
fixed $p,r,M,N$.}

\vskip .1 cm

\paragraph{\bf Proof} Assume that the lemma is not true for some
$0<r<1$, $p>0$ and $N>M>0$. Then there exists $\delta >0$, there
exist two sequences $\varphi_n,\psi_n$ of holomorphic functions
which satisfy $\| \varphi_n\|_{Bers}, \| \psi_n\|_{Bers}\leq N$,
$\|\varphi_n|_{\Delta_r}\|_{Bers}\geq M$ and, there exists a
sequence of measurable set $S_n\subset\Delta_1$, $|S_n|=p$, such
that $\| Arg(\frac{\psi_n}{\varphi_n}
)|_{S_n}\|_{\infty}\leq\frac{1}{n}$ and for each $k>0$ there
exists $z=z(k)\in\Delta_r$ with
\begin{equation}
\label{ineq} |\psi_n(z)-k\varphi_n(z)|\geq \delta >0.
\end{equation}
We find a contradiction with the above statement.

\vskip .3 cm

Since $\|\varphi_n\|_{Bers}, \|\psi_n\|_{Bers}\leq N$, there exist
convergent subsequences $\varphi_{n_k}\to\varphi$,
$\psi_{n_k}\to\psi$ with $\|\varphi \|_{Bers} ,\|\psi \|_{Bers}
<\infty$. The convergence is uniform on compact subsets of
$\Delta_1$. For simplicity of notation write $\varphi_n$,$\psi_n$
in place of $\varphi_{n_k},\psi_{n_k}$. There exists $r_1$,
$r<r_1<1$, such that $|\Delta_{r_1}\cap S_n|\geq\frac{p}{2}$ for
all $n$. Also $\varphi_n\to\varphi$ and $\psi_n\to\psi$ uniformly
on $\Delta_{r_1}$, namely $\| (\varphi_n-\varphi
)|_{\Delta_{r_1}}\|_{Bers},\| (\psi_n -\psi
)|_{\Delta_{r_1}}\|_{Bers}\to 0$ as $n\to\infty$.

\vskip .3 cm

Since $\|\varphi_n|_{\Delta_r}\|_{Bers}\geq M$ for each $n$, we
have $\|\varphi|_{\Delta_r}\|_{Bers}\geq M$ and in particular
$\varphi$ is not identically equal to zero. If $\psi\equiv 0$,
then the above inequality (\ref{ineq}) fails by taking $k>0$ small
enough and $n$ large enough.

\vskip .3 cm

Assume that $\psi$ is not a zero function. The number of zeros of
$\varphi$ and $\psi$ in $\Delta_{r_1}$ is finite. Let $R$ be the
union of disk neighborhoods of the zeros small enough such that
$|S_n\cap (\Delta_{r_1}-R)|\geq\frac{p}{3}$. Given $q>0$, define
$$
D_{q}=\{ z\in\Delta_{r_1}-R;\ |Arg(\frac{\psi}{\varphi})(z)|\leq
q\}
$$
and
$$
D_{q}^n=\{ z\in\Delta_{r_1}-R;\
|Arg(\frac{\psi_n}{\varphi_n})(z)|\leq q\},
$$
and let $D_0=\cap_{k=1}^{\infty}D_{\frac{1}{k}}$ and
$D_0^n=\cap_{k=1}^{\infty}D_{\frac{1}{k}}^n$.

\vskip .3 cm

There exists $n_0$ such that $\| (\varphi_n-\varphi
)|_{\Delta_{r_1}}\|_{Bers}\leq q$ and $\| (\psi_n -\psi
)|_{\Delta_{r_1}}\|_{Bers}\leq q$ for all $n>n_0$. Then there
exists a universal constant $c>0$ such that $D_{cq}\supset
D_{q}^n$ for $n>n_0$. In addition, $D_q^n\supset
S_n\cap(\Delta_{r_1}- R)$ for $n>n_0$ whenever $\frac{1}{n_0}< q$.
Therefore $|D_{cq}|\geq\frac{p}{3}$, for each $q>0$. By the
monotonicity of a positive measure, we obtain
$|D_0|=\lim_{n\to\infty} |D_\frac{1}{n}|\geq \frac{p}{3}$.

\vskip .3 cm

We have $Arg(\frac{\psi}{\varphi})(z)=0$ for all $z\in D_0$. Since
$|D_0|>0$, we obtain that $\varphi (z)=k\psi (z)$ for a fixed
$k>0$ and for all $z\in\Delta_1$. But then $\|
(\varphi_n-k\psi_n)|_{\Delta_{r_1}}\|_{Bers}\to 0$ as
$n\to\infty$, which again gives a contradiction with (\ref{ineq}).
$\Box$

\vskip .3 cm

Let $A_1$ be the unit sphere in $L^1$-norm of the space $A(\H )$
of transversely continuous holomorphic quadratic differentials on
the solenoid $\H$. For a given natural number $N$, define
$$
A_1(N)=\{\varphi\in A_1;\ \|\varphi\|_{Bers}\leq N\}.
$$

\vskip .3 cm

Let $X \subset T(\H)$ be the set of points that do have a Teichm\"uller representative.
For $[\mu ]\in X$ let $k\frac{|\varphi|}{\varphi}$,
$0<k<1$, $\varphi\in A_1$, be that representative. Define the map $\pi :X\to A_1, $ by
$$
\pi ([\mu ])=\varphi,
$$
for $[\mu ]\neq [0]$ and
$$
\pi ([0])=0.
$$

\vskip .3 cm

We need the following proposition.
\vskip .3 cm

\paragraph{\bf Proposition 4.2} {\it Let $[\mu] \in X$, where $\mu$ is a Beltrami
coefficient on $\H$. If $[\mu]$ is an element to the closure of
the set $\pi^{-1} (A_1(N))\cup\{ [0]\}$ then $[\mu] \in  \pi^{-1}
(A_1(N))\cup\{ [0]\}$}.

\vskip .3 cm

\paragraph{\bf Proof}
Let $[\mu_n]\in\pi^{-1}(A_1(N))$ such that $[\mu_n]\to [\mu ]$ in
the Teichm\"uller metric. We need to show that $[\mu ]\in
\pi^{-1}(A_1(N))\cup\{ [0]\}$. If $[\mu ]=[0]$ then we are done.
Therefore, we assume that $[\mu ]\neq [0]$. Without loss of
generality, we assume that
$\mu_n=k_n\frac{|\varphi_n|}{\varphi_n}$ and $\mu =k\frac{|\varphi
|}{\varphi}$ for $\varphi_n,\varphi\in A_1$. Then $k_n\to k$ as
$n\to\infty$ by our assumption. We have that $\varphi_n \in
A_1(N)$ and we need to show that $\varphi\in A_1(N)$.

\vskip .3 cm

There exist $\nu_n\in [\mu_n ]$ such that $\|
Belt((f^{\nu_n})^{-1}\circ f^{\mu})\|_{\infty}\to 0$ as
$n\to\infty$ because $[\mu_n]\to [\mu ]$. Then $\nu_n =\mu +o(1)$,
where $\| o(1)\|_{\infty}\to 0$ as $n\to\infty$. Also
$\|\nu_n\|_{\infty}>k_n$ by the unique extremality of $\mu_n$
\cite{Sa}, and we apply the $\delta$-inequality of \cite{BLMM} to
$\nu_n$ and $\mu_n$:
$$
\int_{\S}\Big{|}\frac{\tilde{\nu}_n(f_n)-\tilde{\mu}_n(f_n)}{1-
\overline{\tilde{\nu}_n(f_n)} \tilde{\mu}_n(f_n)}\Big{|}^2 |\psi
|dm\leq C\Big{(}\|\nu_n\|_{\infty}- Re\int_{\S}\nu_n\psi
dm\Big{)},
$$
for all $\psi\in A_1$, where $f_n=f^{\nu_n}$ and where
$\tilde{\nu}_n,\tilde{\mu}_n$ are the Beltrami coefficients of the
inverse maps of $f^{\nu_n}, f^{\mu_n}$, respectively. The proof of
the $\delta$-inequality for the solenoid follows the same lines as
the proof for Riemann surfaces using the Reich-Strebel inequality
for the solenoid \cite{Sa}. Note that
$\tilde{\nu}_n(f_n)=-\nu_n\lambda_{\nu_n}$ and
$\tilde{\mu}_n(f_n)=-\mu_n\lambda_{\nu_n}+o(1)$, where
$\lambda_{\nu_n}=\frac{\overline{(f_n)_z}}{(f_n)_z}$.

\vskip .3 cm

Then we obtain
$$
\int_{\S}|\mu -\mu_n|^2|\psi |dm\leq C_1(k-Re\int_{\S}\mu\psi
dm)+o(1).
$$
We let $\psi =\varphi$ in the above and obtain
$$
\int_{\S}\Big{|}\frac{|\varphi |}{\varphi}
-\frac{|\varphi_n|}{\varphi_n}\Big{|}^2|\varphi |dm\to 0
$$
as $n\to\infty$. This implies that
$$
Re\int_{\S}\frac{|\varphi_n|}{\varphi_n}\varphi dm\to 1
$$
as $n\to\infty$.

\vskip .3 cm

By the definition $\varphi =\pi ([\mu ])$. By our assumption,
$\|\varphi_n\|_{Bers}\leq N$ and it is enough to show that
$\|\varphi\|_{Bers}\leq N$. Denote by $\alpha $ the product of
leafwise hyperbolic area measure and the transverse measure on
$\S$. We scale the transverse measure in such fashion that $\alpha
(\S )=1$. Define $\S_{n,\epsilon} =\{ x\in\S
:|Arg(\frac{\varphi_n(x)}{\varphi (x)})|<\epsilon\}$. We show that
$\lim_{n\to\infty}\alpha (\S_{n,\epsilon})=1$, for all $\epsilon
>0$. Assume on the contrary that
$\lim\inf_{n\to\infty}\alpha (\S_{n,\epsilon})\leq 1-\delta$, for
$\delta >0$. Then we have
$$
Re\int_{\S}\frac{|\varphi_n |}{\varphi_n}\varphi
dm\leq\int_{\S_{n,\epsilon}}|\varphi |dm+\cos\epsilon\int_{\S
-\S_{n,\epsilon}}|\varphi|dm.
$$

\vskip .3 cm

The zeros of $\varphi$ make a closed subset of $\S$ which is
leafwise discrete and whose $\alpha$-measure is $0$. Moreover,
there exists an open neighborhood $U$ of the zeros of arbitrary
small $\alpha$-measure whose intersection with any leaf consists
of hyperbolic disks. For any such neighborhood $U$, we have
$C(U)=\inf_{(z,t)\in \S -U} |\varphi(z,t)\rho^{-2}(z,t)|>0$. From
the above we get, for $U$ small enough, that
$$
\lim\inf_{n\to\infty}Re\int_{\S}\frac{|\varphi_n
|}{\varphi_n}\varphi dm\leq 1-(1-\cos\epsilon )(\delta -\alpha
(U))C(U)<1
$$
which contradicts $\int_{\S}\frac{|\varphi_n |}{\varphi_n}\varphi
dm\to 1$ as $n\to\infty$. Therefore $\lim_{n\to\infty}\alpha
(\S_{n,\epsilon})=1$ as $n\to\infty$.

\vskip .3 cm

We fix $\delta >0$. By the above, there exists $n_1$ such that
$\alpha (\S_{n_1,\delta /2})>1-\delta /2$. Further, there exists
$n_2\geq n_1$ such that $\alpha (\S_{n_2,\delta /4})>1-\delta /4$,
and so on. In general, we find $n_j$ such that $n_j\geq n_{j-1}$
and $\alpha (\S_{n_j ,\delta /2^j})>1-\delta /2^j$. Define
$\S_0=\cap_{j=1}^{\infty}\S_{n_j,\delta /2^j}$. Then $\alpha (\S_0
) \geq 1-\delta$ and the sequence
$Arg(\frac{\varphi_{n_j}}{\varphi})$ converges to zero uniformly
on $\S_0$. For simplicity of notation, rename the sequence
$\varphi_{n_j}$ to $\varphi_n$.

\vskip .3 cm

Consider lifts $\tilde{\varphi },\tilde{\varphi}_n$ of $\varphi
,\varphi_n$ to the universal covering $\Delta_1\times\hat{G}$ of
$\S$. Recall that for a finite index subgroup $G_k$ of $G$ we have
$\Delta_1\times\hat{G}_k /G_k\equiv\S$. Given $\epsilon >0$ there
exists $k$ such that
$\sup_{z\in\Delta_1}|\tilde{\varphi}(z,t_1)-\tilde{\varphi
}(z,t_2)|\rho^{-2}(z)<\epsilon$, for all $t_1,t_2\in\hat{G}_k$.
Let $\omega_k$ be a fundamental polygon for $G_k$. Then
$\omega_k\times\hat{G}_k$ is a fundamental set for the action of
$G_k$ on $\Delta_1\times\hat{G}_k$.

\vskip .3 cm

Let $c_k=\mbox{h-area}(\omega_k )$. Then $m(\hat{G}_k)=1/c_k$
because of the normalization $\alpha (\S )=1$. Given
$t\in\hat{G}_k$, denote by $\omega_{k,t}^{\delta ,n}$ the set of
all $(z, t)\in\omega_k\times\hat{G}_k$ such that
$|Arg(\frac{\tilde{\varphi}_n}{\tilde{\varphi}})(z,t)|<\delta$.
Let $\hat{G}_k^{\delta_1,n} =\{ t\in\hat{G}_k;\ \mbox{h-area}
(\omega_{k,t}^{\delta ,n})< \delta_1\}$. Then we obtain
$$
1-\delta\leq\alpha (\S_0 )\leq \delta_1
m(\hat{G}_k^{\delta_1,n})+c_k(1/c_k-m(\hat{G}_n^{\delta_1,n}))
$$
which implies
$$
m(\hat{G}_k^{\delta_1,n})\leq\frac{\delta}{c_k-\delta_1}.
$$

\vskip .3 cm

The above implies that $\omega_{k,t}^{\delta ,n}$ has Lebesgue
measure bounded from below for each $t\in
\hat{G}_k-\hat{G}_k^{\delta_1,n}$ and Lemma 4.1 applies to such
$t$. Thus, by Lemma 4.1, there exists a sequence of functions
$k_n:\hat{G}_k-\hat{G}_k^{\delta_1,n}\to\mathbb{R}^{+}$ such that
\begin{equation}
\label{lim} \|
(\tilde{\varphi}_n-k_n\tilde{\varphi})|_{\omega_k\times
(\hat{G}_k-\hat{G}_k^{\delta_1,n})}\|_{Bers}\leq d_n\to 0
\end{equation}
 as $n\to\infty$.

\vskip .3 cm

We claim that there exists a sequence
$t_n\in\hat{G}_k-\hat{G}_k^{\delta_1,n}$ such that
$\lim\sup_{n\to\infty}k_n(t_n)\geq 1$. Suppose on the contrary
that there exists $c>0$ such that $\| k_n\|_{\infty}\leq 1-c$ for
all large $n$. From (\ref{lim}) and by the above assumption, we
obtain
$$
\int_{\omega_k\times
(\hat{G}_k-\hat{G}_k^{\delta_1,n})}|\tilde{\varphi}_n|dm-d_n\leq
\int_{\omega_k\times
(\hat{G}_k-\hat{G}_k^{\delta_1,n})}k_n|\tilde{\varphi}|dm\leq
(1-c)\int_{\omega_k\times
(\hat{G}_k-\hat{G}_k^{\delta_1,n})}|\tilde{\varphi}|dm.
$$
This implies
\begin{equation}
\label{contr}
 \int_{\omega_k\times
(\hat{G}_k-\hat{G}_k^{\delta_1,n})}|\tilde{\varphi}_n|dm/(1-c)-d_n/(1-c)
\leq \int_{\omega_k\times
(\hat{G}_k-\hat{G}_k^{\delta_1,n})}|\tilde{\varphi}|dm.
\end{equation}

\vskip .3 cm

Note that
\begin{equation}
\label{contr1}
\int_{\omega_k\times\hat{G}_k^{\delta_1,n}}|\tilde{\varphi}_n|dm
\leq \|\varphi_n\|_{Bers}\ \alpha
(\omega_k\times\hat{G}_k^{\delta_1,n})\leq
Nc_k\frac{\delta}{c_k-\delta_1}\to 0
\end{equation}
as $\delta\to 0$, for fixed $k$ and $\delta_1$, and uniformly in
$n$. If we take $n$ large enough and $\delta$ small enough in
(\ref{contr}) and (\ref{contr1}), we get that $\int_{\S}|\varphi
|dm>1$ which is a contradiction with our choice of $\varphi$.

\vskip .3 cm

Therefore, there exists $t_n$ such that
$\lim\sup_{n\to\infty}k_n(t_n)\geq 1$. From (\ref{lim}) we get
that $\|\tilde{\varphi}|_{\omega_k\times t_n}\|_{Bers}\leq
\|\tilde{\varphi}_n|_{\omega_k\times
t_n}\|_{Bers}/k_n(t_n)+d_n/k_n(t_n).$ Consequently, we have that
$\|\varphi\|_{Bers}\leq N+\epsilon$ by letting $n\to\infty$. Since
$\epsilon$ was arbitrary, we get that $\|\varphi\|_{Bers}\leq N$.
$\Box$

\vskip .3 cm

Before proving the Theorem 1. we need to prove the next lemma.
\vskip .3cm
Let $S_0$ be a compact Riemann surface of genus two at
least two. Let $\gamma$ be a non-diving simple closed geodesic (in
the corresponding hyperbolic metric). We cut $S_0$ along $\gamma$
to obtain a bordered hyperbolic surface $S_0^b$. Denote by $S_n$ a
compact Riemann surface obtained by gluing $n$ copies of $S_0^b$
along their boundaries such that $S_n$ is $\mathbb{Z}_n$-cover of
$S_0$. Given $r$, $0<r<1$, we denote by $R_{n,r}$ a subsurface of
$S_n$ which consists of $[rn]$ consecutive copies of $S_0^b$ in
$S_n$, where $[rn]$ is the smallest integer less that or equal to
$rn$. The hyperbolic metric on $S_0$ lifts to a unique hyperbolic
metric on $S_n$ and both hyperbolic metrics lift to a transversely
locally constant hyperbolic metric on $\S$. Denote by
$\tilde{R}_{n,r}$ the lift of $R_{n,r}$ to $\S$. It is then clear
that $\alpha (\tilde{R}_{n,r} )=[nr]/n$.

\vskip .3 cm

If $f$ is a leafwise quadratic differential on $\S$, we define $\|
f\|_{L^1}:=\int_{\S}|f|dm$. Define a non-holomorphic quadratic
differential on $\S$ by
\begin{equation*}
\tilde{\varphi}_n= \left\{
\begin{array}l
\ 0,\ \ \mbox{ on   }\ \ \S -\tilde{R}_{n,r}\\
\tilde{\varphi}_0,\ \mbox{ on   }\ \ \tilde{R}_{n,r}
\end{array}
\right.
\end{equation*}
Note that $\tilde{\varphi}_n$ is the lift of
\begin{equation*}
\varphi_n= \left\{
\begin{array}l
\ 0,\ \ \mbox{ on   }\ \ S_n -R_{n,r}\\
\varphi_0 ,\ \mbox{ on   }\ \ R_{n,r}
\end{array}
\right.
\end{equation*}

\vskip .3 cm

\paragraph{\bf Lemma 4.3}{\it There exists a holomorphic quadratic
differential $\tilde{\psi}_n\in A(\S )$ such that
$$
\Big{\|}\frac{\tilde{\psi}_n}{\|\tilde{\psi}_n\|_{L^1}}
-\frac{\tilde{\varphi}_n}{\|\tilde{\varphi}_n\|_{L^1}}\Big{\|}_{L^1}\to
0,
$$
as $n\to\infty$, for fixed $r$. Moreover, we can choose
$\tilde{\psi}_n$ to be the lift of a holomorphic quadratic
differential $\psi_n$ on $S_n$.}

\vskip .3 cm

\paragraph{\bf Proof}
Let $A(S_n)$ denote the space of holomorphic quadratic
differentials on $S_n$. We define a linear functional $\sigma :
A(S_n)\to \mathbb{C}$ by
$$
\sigma (f )=\int_{S_n}\bar{\varphi}_n\rho^{-2}f,
$$
for $f\in A(S_n)$. It is a standard fact for Riemann surfaces that
there exists a unique $\psi_n\in A(S_n)$ such that $\sigma
(f)=\int_{S_n}\bar{\psi}_n\rho^{-2}f$, for all $f\in A(S_n)$.
Denote by $\tilde{\psi}_n$ the lift of $\psi_n$ to $\S$.

We consider $c$-neighborhood $U_n(c)$ of the two boundary curves
of $R_{n,r}$ in the hyperbolic metric on $S_n$, for $c>0$. We
claim that for any $\epsilon >0$ there exist $n_0$ and $c>0$ such
that
\begin{equation}
\label{un} \rho^{-2}|\varphi_n-\psi_n|<\epsilon
\end{equation}
on $S_n-U_n(c)$ for all $n>n_0$.

To show the claim, we assume that it is not true (and arrive at a
contradiction). Then there exist $\epsilon >0$, a sequence $c_n>0$
and a sequence $z_n\in S_n-U_n(c_n)$ such that $c_n\to\infty$ as
$n\to\infty$ and
\begin{equation}
\label{three} \rho^{-2}(z_n)|\varphi_n
(z_n)-\psi_n(z_n)|\geq\epsilon
\end{equation}
for all $n$. We arrange that either $z_n\in S_n-R_{n,r}$ or
$z_n\in R_{n,r}$ for all $n$ after possibly taking a subsequence.
Consider a $\mathbb{Z}$-cover $\tilde{S}$ of $S_0$ which is made
by gluing together infinitely many $S_0^b$. Note that $\tilde{S}$
is also $\mathbb{Z}$-cover of each $S_n$. We arrange that the
covering maps $\tilde{S}\to S_n$ have a lift of $z_n$ in a fixed
copy of $S_0^b$ in $\tilde{S}$. We denote by $\tilde{\varphi}_n$,
$\tilde{\psi}_n$ the lifts of $\varphi_n$, $\psi_n$ to $\tilde{S}$
as well as to $\S$ and the meaning should be read from the
context. Define a linear functional on $A(\tilde{S} )$ by
$$\tilde{\sigma}_n(f):=\int_{\tilde{S}}\stackrel{-} {\tilde{\varphi}}_n\rho^{-2}f,$$
for all $f\in A(\tilde{S} )$. Then $\tilde{\sigma}_n$ satisfies
$$\tilde{\sigma}_n(f)=\int_{\tilde{S}}\stackrel{-}
{\tilde{\psi}}_n\rho^{-2}f$$ for all $f\in A(\tilde{S })$ by the
Bers' reproducing formula. Note that $A(S_n)$ does not lift to a
subset of $A(\tilde{S})$.

By the choice of the covering maps $\tilde{S}\to S_n$, we have
that either $\tilde{\varphi}_n\to\tilde{\varphi}_0$ on $\tilde{S}$
uniformly on compact subsets if $z_n\in R_{n,r}$, or
$\tilde{\varphi}_n\to 0$ uniformly on compact subsets otherwise.
Then either
$$
\lim_{n\to\infty}\int_{\tilde{S}}\stackrel{-}
{\tilde{\varphi}}_n\rho^{-2} f=\int_{\tilde{S}} \stackrel{-}
{\tilde{\varphi}}_0\rho^{-2} f
$$
for all $f\in A(\tilde{S})$ in the first case or
$$
\lim_{n\to\infty}\int_{\tilde{S}}\stackrel{-}
{\tilde{\varphi}}_n\rho^{-2} f=0
$$
in the second case.

It is clear that $\|\tilde{\psi}_n\|_{Bers}\leq
3\|\varphi_0\|_{Bers}<\infty$ for all $n$ because
$\|\tilde{\sigma}_n\|\leq \|\varphi_0\|_{Bers}$. Therefore,
$\tilde{\psi}_n$ has a subsequence which converges uniformly on
compact subsets of $\tilde{S}$ to a holomorphic quadratic
differential $\tilde{\psi}\in A(\tilde{S})$. Then
\begin{equation}
\label{four} \lim_{n\to\infty}\int_{\tilde{S}}\stackrel{-}
{\tilde{\psi}}_n\rho^{-2} f=\int_{\tilde{S}} \stackrel{-}
{\tilde{\psi}}\rho^{-2} f.
\end{equation}
On the other hand, $\tilde{\psi}\neq\tilde{\varphi}_0$ in the
first case and $\tilde{\psi}\neq 0$ in the second case by
(\ref{three}) and the fact that the inequality prevails in the
lifts to a compact subset of $\tilde{S}$. Thus we obtain two
different presentations for the limiting linear functional. This
is a contradiction to the uniqueness of the presentation of linear
functionals in the above form. Therefore, we showed that given
$\epsilon >0$ there exists $c>0$ such that $\rho^{-2}|\varphi_n-
\psi_n|<\epsilon$ in $S_n-U_n(c)$ for all $n$ large enough.

Note that the hyperbolic area of $U_n(c)$ is constant in $n$ for a
fixed $c>0$ by our choice of coverings. Since the genus of $S_n$
goes to infinity as $n\to\infty$, we conclude that $\alpha
(\tilde{U}_n(c))\to 0$ as $n\to\infty$, where $\tilde{U}_n(c)$ is
the lift of $U_n(c)$ to $\S$. Then by (\ref{un}) and by the above,
we get
\begin{equation*}
\begin{array}l
\int_{\S}|\tilde{\varphi}_n-\tilde{\psi}_n|dm\leq \int_{\S
-(\tilde{U}_n(c))}|\tilde{\varphi}_n-\tilde{\psi}_n|\rho^{-2}\rho^2dm
+\int_{\tilde{U}_n(c)}|\tilde{\varphi}_n-\tilde{\psi}_n|\rho^{-2}\rho^2dm
\leq
\\ \ \ \ \leq
\epsilon\alpha (\S -U_n(c))+M\alpha (U_n(c))\to 0
\end{array}
\end{equation*}
as $n\to\infty$. In other words, we showed that $\lim_{n\to\infty}
\|\tilde{\varphi}_n-\tilde{\psi}_n\|_{L^1}=0$.

It is clear that
$\|\tilde{\varphi}_n\|_{L^1}=\frac{[rn]}{n}\|\varphi_0\|_{L^1(S_0)}$.
By the above, we also get that $\|\tilde{\psi}_n\|_{L^1}-
\frac{[rn]}{n}\|\varphi_0\|_{L^1(S_0)}\to 0$ as $n\to\infty$. This
implies that both $\|\tilde{\varphi}_n\|_{L^1}$ and
$\|\tilde{\psi}_n\|_{L^1}$ are bounded from below independently of
$n$. Moreover, their difference converges to $0$ as $n\to\infty$.
Thus we obtain
\begin{equation*}
\Big{\|}\frac{\tilde{\psi}_n}{\|\tilde{\psi}_n\|_{L^1}}-
\frac{\tilde{\varphi}_n}{\|\tilde{\varphi}_n\|_{L^1}}\Big{\|}_{L^1}\leq
\frac{k[nr]}{n\|\varphi_0\|_{L^1(S_0)}}\|\tilde{\varphi}_n-\tilde{\psi}_n\|_{L^1}
\to 0,
\end{equation*}
for some constant $k>0$ as $n\to\infty$. $\Box$

\vskip .5 cm

\paragraph{\bf Proof of Theorem 1}

To prove that the set of points which do not have a Teichm\"uller representative is generic we need to prove that
the set $X$ is of first category. Recall that

\begin{equation*}
X=\{ [0]\}\cup (\cup\pi^{-1}(A_1(N))),
\end{equation*}
where the second union is over all ${N\in\mathbb{N}}$. It is enough to show that each of the sets
$\pi^{-1}(A_1(N))$ is of the first category. We prove this by contradiction.
\vskip .3cm
Assume that for some $N$ the set $\pi^{-1}(A_1(N))$ is of the second category.
Then the closure $\pi^{-1}(A_1(N))^c$ has non-empty interior. Moreover, by Proposition 4.2 every
element in $\pi^{-1}(A_1(N))^c$ that has a Teichm\"uller representative must be in
$\pi^{-1}(A_1(N))$. Let $[\mu=k\frac{|\varphi |}{\varphi}]$
be a transversely locally constant Beltrami coefficient on $\S$ such that the point $[\mu]$
is an element of the interior of the set $\pi^{-1}(A_1(N))^c$.
We assume that $\mu$ is lifted from a surface $S_0$, that is $\varphi\in A_1$ is a lift
of a holomorphic quadratic differential $\varphi_0$ on a closed Riemann surface $S_0$ of genus at least two.
We will show that there exists a transversely locally constant sequence
$[\mu_n]\to [\mu ]$ such that $\xi_n=\pi ([\mu_n])$ are unbounded
in the Bers norm. Since for $n$ large enough we have that $[\mu_n]$ is an element of the interior of the set
$\pi^{-1}(A_1(N))^c$ we will  obtain a contradiction.
\vskip .3cm
We keep the notation $S_n$ for $\mathbb{Z}_n$ cover of
$S_0$ and $R_{n,r}$ for $[rn]/n$ proportion of $S_n$ as above. Consider a Beltrami coefficient
\begin{equation}
\label{mun} \mu_n=\Big{\{}
\begin{array}l
\ \ \ \ k\frac{|\varphi |}{\varphi},\ \ \ \ \ \ \ \ \S -\tilde{R}_{n,r}\\
(1+r)k\frac{|\varphi |}{\varphi},\ \ \ \ \tilde{R}_{n,r}
\end{array}
\end{equation}
for $r>0$ small enough such that $\|\mu_n\|_{\infty}<1$. The
Beltrami coefficient $\mu_n$ is not smooth at the lift of two
boundary curves of $R_{n,r}$ to the solenoid $\S$ and it can be
smoothly approximated in arbitrary small area neighborhoods of the
lift. Therefore, we can work with $\mu_n$ as well and we refer the
reader to \cite{Sa} for more details.

\vskip .3 cm

Let $\lambda^n(\psi )=\int_{\S}\mu_n\psi dm$, for $\psi\in A(\S
)$, be the corresponding linear functional. Denote by
$\tilde{\varphi}_n$ the non-holomorphic quadratic differential
which is the lift of
\begin{equation}
\label{hq} \varphi_n=\Big{\{}
\begin{array}l
\varphi_0,\ \ R_{n,r}\\
0,\ \ \ \ S_n-R_{n,r}
\end{array}.
\end{equation}

\vskip .3 cm

By Lemma 4.3, there exist holomorphic quadratic differentials
$\tilde{\psi}_n$ on $\S$ such that
\begin{equation}
\label{l1con}
\Big{\|}\frac{\tilde{\psi}_n}{\|\tilde{\psi}_n\|_{L^1}}-
\frac{\tilde{\varphi}_n}{\|\tilde{\varphi}_n\|_{L^1}}\Big{\|}_{L^1}
\to 0
\end{equation}
as $n\to\infty$.

\vskip .3 cm

Denote by $\|\lambda_n\| =\sup_{\psi\in A_1} |\lambda_n(\psi )|$
the operator norm of $\lambda_n$. Then $\|\lambda_n\|\leq (1+r)k$
because $\|\mu_n\|_{\infty} =(1+r)k$. Since
$\lambda_n(\tilde{\varphi}_n/\|\tilde{\varphi}_n\|_{L^1})=(1+r)k$
and by (\ref{l1con}), we have that $\|\lambda_n\|\to (1+r)k$ as
$n\to\infty$. There exists $\xi_n\in A_1$ (which is the lift of a
holomorphic quadratic differential on $S_n$) and there exists
$l_n>0$ such that $l_n\frac{|\xi_n|}{\xi_n}\in [\mu_n]$.

\vskip .3 cm

Let $k(\mu_n )=\inf_{\nu\in [\mu_n]}\|\nu\|_{\infty}$. The
Teichm\"uller contraction inequality \cite{Sa} applied to $\mu_n$
gives
$$
\|\mu_n\|_{\infty}-k(\mu_n)\leq C(\|\mu_n\|_{\infty}
-\sup_{\psi\in A_1}Re\int_{\S}\mu_n\psi dm)
$$
where $C>0$ is a fixed constant. The right hand side of the above
inequality converges to zero as $n\to\infty$. Therefore the left
hand side converges to zero as well. Since
$\|\mu_n\|_{\infty}=(1+r)k$, then we have $k(\mu_n)=l_n\to (1+r)k$
as $n\to\infty$.

\vskip .3 cm

We use another standard formula which is an easy consequence of
the Reich-Strebel inequality developed for the solenoid in
\cite{Sa} in the course of proof of the Teichm\"uller contraction.
Namely, we get
$$
\frac{1+l_n}{1-l_n}\leq\int_{\S}\frac{|1-\mu_n\frac{\xi_n}{|\xi_n
|}|^2}{1-|\mu_n|^2}|\xi_n|dm
$$
where $l_n\frac{|\xi_n|}{\xi_n}\in [\mu_n]$ is an important
condition (the formula is not true for arbitrary holomorphic
quadratic differential). Further,
$$
\int_{\S}\frac{|1-\mu_n\frac{\xi_n}{|\xi_n
|}|^2}{1-|\mu_n|^2}|\xi_n|dm \leq \int_{\S
-\tilde{R}_{n,r}}\frac{(1+k)^2}{1-k^2}|\xi_n|dm
+\int_{\tilde{R}_{n,r}}\frac{[1+(1+r)k]^2}{1-[(1+r)k]^2}|\xi_n|dm
$$
which implies that
$$
\frac{1+l_n}{1-l_n}\leq \frac{1+k}{1-k}\int_{\S
-\tilde{R}_{n,r}}|\xi_n|dm+
\frac{1+(1+r)k}{1-(1+r)k}\int_{\tilde{R}_{n,r}} |\xi_n|dm.
$$

\vskip .3 cm

From the above inequality and by $l_n\to (1+r)k$, we get that
$\lim_{n\to\infty}\int_{\S -\tilde{R}_{n,r}}|\xi_n|dm=0$, for each
$0<r<1$. To see this assume on the contrary that
${\lim\sup}_{n\to\infty}\int_{\S -\tilde{R}_{n,r}}|\xi_n|dm=\delta
>0$. Then by taking $\lim\sup_{n\to\infty}$ in the above
inequality, we obtain $
\frac{1+l_n}{1-l_n}\leq\frac{1+k}{1-k}\delta
+\frac{1+(1+r)k}{1-(1+r)k}(1-\delta )$ which is impossible. By
Cantor diagonal argument, there exists a sequence $r_n\to 0$,
$0<r_n<1$, such that $nr_n\to \infty$ and $\int_{\S
-\tilde{R}_{n,r_n}}|\xi_n|dm\to 0$ as $n\to\infty$. Clearly
$\lim_{n\to\infty}\alpha (\tilde{R}_{n,r_n})=0$. For simplicity,
we write $\tilde{R}_n=\tilde{R}_{n,r_n}$.

\vskip .3 cm

Finally, we assume that $\|\xi_n\|_{Bers}\leq N$. Then we have
that
$$
1=\|\xi_n\|_{L^1}\leq \int_{\S
-\tilde{R}_n}|\xi_n|dm+\|\xi_n\|_{Bers}\ \alpha
(\tilde{R}_n)=\int_{\S -\tilde{R}_n}|\xi_n|dm+N \alpha
(\tilde{R}_n) \to 0
$$
as $n\to\infty$. This is a contradiction. Therefore $\xi_n$ is not
bounded in the Bers norm.

$\Box$

\section{Infinitesimal Teichm\"uller classes without Teichm\"uller
representatives}

A tangent vector to $T (\S )$ at the basepoint $[0]$ is
represented by a Beltrami differential $\mu$. It defines a
continuous linear functional on $A(\S )$ via the natural pairing.
If $\mu$ is infinitesimally equivalent to $k\frac{|\varphi
|}{\varphi}$ (i.e. $\mu -k\frac{|\varphi |}{\varphi}\in N(\S )$)
then the linear functional achieves norm on the unique vector
$\frac{\varphi}{\|\varphi\|_{L^1}}\in A(\S )$. A question whether
each infinitesimal Teichm\"uller class contains a Beltrami
coefficient of the Teichm\"uller type $k\frac{|\varphi
|}{\varphi}$, for $k\in\mathbb{R}^{+}$ and $0\neq\varphi\in A(\S
)$, is analogous to the question of the existence of Teichm\"uller
representative for marked complex structures. An equivalent
question is whether the induced linear functional achieves its
norm on $A(\S )$. We obtain an analogous result to Theorem 1:

\vskip .5 cm

\paragraph{\bf Theorem 3} {\it The set of points in the tangent
space of $T(\S )$ at the basepoint $[0]$ which do not achieve its norm
on $A(\S )$ is generic.}

\vskip .3 cm

\paragraph{\bf Proof} Denote by $B(\S )=L^{\infty}(\S )/N(\S)$
the tangent space to $T(\S )$ at the basepoint. Let $X \subset B(\S )$ be the set of points
that do achieve its norm. That is, each $\lambda\in X$ achieves its norm on some
$\varphi\in A_1$. We define $\pi :X\to A_1$ by $\pi (\lambda
)=\varphi$ if $\lambda$ achieves its norm on $\varphi\in A_1$.
Then we have $X=\cup_{N=1}^{\infty}\pi^{-1}(A_1(N))\cup\{[0]\} )$.
\vskip .3cm
We show that if $\lambda \in {\pi^{-1}(A_1(N))} ^{c}$ achieves its norm then $\lambda \in  \pi^{-1}(A_1(N))$
Let $\lambda_n\in \pi^{-1}(A_1(N))$ and assume that
$\lambda_n \to \lambda$. Let
$\mu_n=k_n\frac{|\varphi_n|}{\varphi_n}$ be the Teichm\"uller
Beltrami differential representing $\lambda_n$, where
$\varphi_n\in A_1$ and let $\mu =k\frac{|\varphi |}{\varphi}$ be
the Teichm\"uller Beltrami differential representing $\mu$, where
$\varphi\in A_1$. Then $k_n\to k$ and
$$
\int_{\S}k_n\frac{|\varphi_n|}{\varphi_n}\psi
dm\to\int_{\S}k\frac{|\varphi |}{\varphi}\psi dm
$$
as $n\to\infty$ for all $\psi\in A(\S )$. By letting $\psi
=\varphi$ in the above, we get
$$
\int_{\S}\frac{|\varphi_n|}{\varphi_n}\varphi dm\to 1
$$
as $n\to\infty$. In the proof of Proposition 4.2, we showed that
the above convergence implies that $\varphi\in A_1(N)$, which proves the claim.

\vskip .3 cm

Same as in the proof of Theorem 1 , we prove that $X$ is of the first category. It is enough to
prove that each $\pi^{-1}(A_1(N))$ is of the first kind.  We do this by contradiction.
\vskip .3cm
Assume that for some $N$ the set $\pi^{-1}(A_1(N))$ is of the second kind. Therefore, the closure
${\pi^{-1}(A_1(N))}^{c}$ has non-empty interior.
Then there exists $\lambda_0 \in (\pi^{-1}(A_1(N)))^{\circ}$ which is transversely locally
constant. This implies that any sequence of transversely locally
constant $\lambda_n$ which converge to $\lambda_0$ must have
bounded Bers norm. We find a contradiction with this statement by
constructing a convergent sequence with unbounded Bers norm below.

\vskip .3 cm

Assume that $\lambda_0\in B(\S )$ achieves its norm on
$\tilde{\varphi}_0\in A_1$, where $\tilde{\varphi}_0$ is the lift
of a holomorphic quadratic differential $\varphi_0$ on a Riemann
surface $S_0$. As in Section 4, define a non-holomorphic quadratic
differential on $\S$ by
\begin{equation*}
\tilde{\varphi}_n= \left\{
\begin{array}l
\ 0,\ \ \mbox{ on   }\ \ \S -\tilde{R}_{n,r}\\
\tilde{\varphi}_0,\ \mbox{ on   }\ \ \tilde{R}_{n,r}
\end{array}
\right.
\end{equation*}
Note that $\tilde{\varphi}_n$ is the lift of
\begin{equation*}
\varphi_n= \left\{
\begin{array}l
\ 0,\ \ \mbox{ on   }\ \ S_n -R_{n,r}\\
\varphi_0 ,\ \mbox{ on   }\ \ R_{n,r}
\end{array}
\right.
\end{equation*}

\vskip .3 cm

By Lemma 4.3, there exist a sequence $\tilde{\psi}_n\in A(\S )$
such that
$\Big{\|}\frac{\tilde{\psi}_n}{\|\tilde{\psi}_n\|_{L^1}}-
\frac{\tilde{\varphi}_n}{\|\tilde{\varphi}_n\|_{L^1}}\Big{\|}_{L^1}\to
0$ as $n\to\infty$. We define
$\psi_n:=\frac{\tilde{\psi}_n}{\|\tilde{\psi}_n\|_{L^1}}$, i.e.
$\psi_n$ is a positive multiple of $\tilde{\psi}_n$ which belongs
to $A_1$.

\vskip .3 cm

Consider linear functionals
$$
\lambda_n(f)=\int_{\S
-\tilde{R}_{n,r}}\frac{|\varphi_0|}{\varphi_0}fdm+(1+l)
\int_{\tilde{R}_{n,r}}\frac{|\varphi_0|}{\varphi_0}fdm.
$$
for $f\in A(\S )$ and $l>0$. Then we obtain
$$
|\lambda_n(\psi_n )|\geq -\int_{\S
-\tilde{R}_{n,r}}|\psi_n|dm+(1+l)\int_{\tilde{R}_{n,r}}|\psi_n|dm
\to 1+l,
$$
as $n\to\infty$. Therefore, $\|\lambda_n\|\to 1+l$ as
$n\to\infty$.

\vskip .3 cm

Note that $\lambda_n$ descends to a functional on $A(S_{n})$. Thus
there exists a unique $\xi_n\in A(S_{n})$ on which $\lambda_n$
achieves its norm. Lift a positive multiple of $\xi_n$ to a
transversely locally constant holomorphic quadratic differential
$\tilde{\xi}_n$ on $\S$ such that $\|\tilde{\xi}_n\|_{L^1}=1$.
Then we have $\lambda_n (\tilde{\xi}_n)=\|\lambda_n\|\geq
1+l-\epsilon$ for all $n$ large enough depending of $\epsilon$.

\vskip .3 cm

We claim that $\int_{\S -\tilde{R}_{n,r}}|\tilde{\xi}_n|dm\to 0$
as $n\to\infty$. Assume on the contrary that there exists $\delta
>0$ such that $\lim\sup_{k\to\infty} \int_{\S
-\tilde{R}_k}|\tilde{\xi}_k|dm=\delta$. Then we get
$$
\lim\sup_{n\to\infty} \lambda_n(\tilde{\xi}_n)\leq \delta
+(1+l)(1-\delta )=1+l-l\delta<1+l.
$$
But this is in contradiction with $\|\lambda_n\| \to 1+l$ as
$n\to\infty$. Therefore $\int_{\S
-\tilde{R}_{n,r}}|\tilde{\xi}_n|dm\to 0$ as $n\to\infty$.

\vskip .3 cm

By the Cantor diagonal argument, there exists a sequence $r_n\to
0$ such that $\int_{\S -\tilde{R}_{n,r_n}} |\tilde{\xi}_n|dm\to 0$
and $nr_n\to \infty$ as $n\to\infty$. The condition $r_n\to 0$
implies that $\alpha (\tilde{R}_{n,r_n})\to 0$ as $n\to\infty$.
For simplicity of notation, we write
$\tilde{R}_n=\tilde{R}_{n,r_n}$.

\vskip .3 cm

This implies
$$
1=\|\tilde{\xi}_n\|_{L^1}=\int_{\S
-\tilde{R}_{n}}|\tilde{\xi}_n|dm+\int_{\tilde{R}_n}|\tilde{\xi}_n|dm\leq
\int_{\S -\tilde{R}_n}|\tilde{\xi}_n|dm
+\|\tilde{\xi}_n\|_{Bers}\alpha (\tilde{R}_n)\to 0
$$

as $n\to\infty$. This is a contradiction.

$\Box$

\vskip .4 cm

\paragraph{\bf Remark} The results of this paper immediately generalize to the
punctured solenoid introduced in \cite{PS}.

\end{document}